\documentclass[twocolumn,twoside]{IEEEtran}

\usepackage[dvipsnames]{xcolor}
\usepackage{amsmath,amssymb,amsthm,epsfig,color,subfigure,empheq,graphicx,pgfplots}
\usepackage{url,algorithm,wasysym,epstopdf,balance,enumitem}
\usepackage{pgfplots}
\pgfplotsset{compat=1.6} 
\usepackage[bb=boondox]{mathalfa}
\usepackage{dsfont}

\usepackage{circuitikz}
\usetikzlibrary{positioning, shapes}
\usetikzlibrary{decorations.pathreplacing}
\usetikzlibrary{arrows}
\usetikzlibrary{plotmarks}
\usetikzlibrary{patterns}
\usetikzlibrary{spy}
\usetikzlibrary{calc,trees,positioning,arrows,chains,shapes.geometric,%
decorations.pathreplacing,decorations.pathmorphing,shapes,decorations.markings,matrix,shapes.symbols}

\tikzset{
>=stealth',
  punktchain/.style={rectangle, rounded corners, 
    draw=black, very thick,text width=10em, 
    minimum height=3em, text centered, on chain},
  line/.style={draw, thick, <-},
  element/.style={tape,top color=white,bottom color=blue!50!black!60!,
    minimum width=8em,draw=blue!40!black!90, very thick,
    text width=10em, minimum height=3.5em, text centered, on chain},
  every join/.style={->, thick,shorten >=1pt},
  tuborg/.style={decorate},
  tubnode/.style={midway, right=2pt},
}
\usepackage[font=footnotesize,skip=5pt]{caption}

\newcommand{\oprocendsymbol}{\hbox{$\bullet$}}
\newcommand{\oprocend}{\relax\ifmmode\else\unskip\hfill\fi\oprocendsymbol}

\newtheorem{theorem}{Theorem}
\newtheorem{proposition}{Proposition}

\newtheorem{remark}{Remark}

\newcommand{\longthmtitle}[1]{\mbox{}{\textbf{\textit{(#1).}}}}

\makeatletter
\def\subsection{\@startsection{subsection}{2}{\z@}{1ex plus 1ex minus 0ex}%
{0.7ex plus .5ex minus 0ex}{\normalfont\normalsize\itshape}}%
\makeatother

\usepackage[absolute,overlay]{textpos}
\usepackage[hyperindex=true, %
  colorlinks=true,%
  pagebackref=false,%
  plainpages=false,%
  pdfpagelabels,%
  linkcolor=NavyBlue,%
  citecolor=NavyBlue,%
  filecolor=NavyBlue,%
  urlcolor=NavyBlue%
]{hyperref}

\begin{document}

\begin{textblock*}{\textwidth}(15mm,9mm) 
\centering \bf \textcolor{NavyBlue}{To appear in the \emph{IEEE Control Systems Letters} \\\url{https://doi.org/10.1109/LCSYS.2022.3228855}}
\end{textblock*}

\title{Dynamic shaping of the grid frequency response shaping using grid-forming IBRs}

\title{Frequency shaping control for weakly-coupled grid-forming IBRs}

\author{Bala Kameshwar Poolla,\, Yashen Lin,\, Andrey Bernstein,\, Enrique Mallada,\, Dominic Gro\ss{}

\thanks{
B. K. Poolla, Y. Lin, and A. Bernstein are with the National Renewable Energy Laboratory (NREL), Golden CO, USA (Email: {\tt\small \{bpoolla, yashen.lin, andrey.bernstein\}@nrel.gov}). Enrique Mallada is with the Electrical and Computer Engineering at the Johns Hopkins University, Baltimore MA, USA (Email: {\tt\small mallada@jhu.edu}). Dominic Gro\ss{} is with the Department of Electrical and Computer Engineering Department at the University of Wisconsin, Madison WI, USA (Email: {\tt\small dominic.gross@wisc.edu}).

This work was authored by the National Renewable Energy Laboratory, operated by Alliance for Sustainable Energy, LLC, for the U.S. Department of Energy (DOE) under Contract No. DE-AC36-08GO28308. Funding provided by NREL Laboratory Directed Research and Development Program. The views expressed in the article do not necessarily represent the views of the DOE or the U.S. Government. The U.S. Government retains and the publisher, by accepting the article for publication, acknowledges that the U.S. Government retains a nonexclusive, paid-up, irrevocable, worldwide license to publish or reproduce the published form of this work, or allow others to do so, for U.S. Government purposes.}%
}	

\maketitle

\begin{abstract}
We consider the problem of controlling the frequency of low-inertia power systems via inverter-based resources (IBRs) that are weakly connected to the grid. We propose a novel grid-forming control strategy, the so-called frequency shaping control, that aims to shape the frequency response of synchronous generators (SGs) to load perturbations so as to efficiently arrest sudden frequency drops. Our solution relaxes several existing assumptions in the literature and is able to navigate trade-offs between peak power requirements and maximum frequency deviations. Finally, we analyze the robustness to imperfect knowledge of network parameters, while particularly highlighting the importance of accurate estimation of these parameters.
\end{abstract}
\section{Introduction}
Electric power systems are undergoing an unprecedented transformation  towards replacing conventional bulk power generation with renewable generation. A significant proportion of emerging technologies, such as photovoltaics, wind power, and battery energy storage systems, are connected to power systems by means of power electronic inverters. In contrast to conventional synchronous generators (SGs) that provide significant kinetic energy storage (i.e., inertia) and primary frequency control via their turbine-governor system, renewables and power electronic inverters deployed today are controlled to maximize the renewable power generation and jeopardize stability \cite{WEB+15,MDH+18}. Instead, grid-forming (GFM) inverters, that impose a well-defined AC voltage are envisioned to replace SGs as the cornerstone of future power systems \cite{MDH+18,MBP+2019}.

The prevalent GFM control strategies in the literature are droop-control \cite{CDA93}, virtual synchronous machines \cite{DSF2015}, and (dispatchable) virtual oscillator control \cite{JDH+14,GCB+19}. While virtual oscillator control is typically aimed at $100\%$ inverter systems, the literature on low-inertia systems that contain a mix of conventional generation and inverter-based resources (IBRs) largely focuses on droop control and virtual synchronous machines that mimic the speed droop and inertia response of SGs to various levels of fidelity. Although strategies based on machine-emulation are compatible with the legacy power system, they do not leverage the full potential of a system of SGs and IBRs. In particular, SGs provide significant (but slow) frequency control, whereas IBRs can provide a fast frequency control response but often cannot sustain this response over a longer duration due to limited energy storage and limited flexibility of their primary power source.

One alternative to inertia emulation or droop control is to leverage the flexible power output of IBRs to \emph{shape} the frequency response of SGs to load changes~\cite{YCVM-2021,YJ-AB-PV-EM:21}. Inspired by model-matching techniques, \emph{frequency shaping control} (FSC) aims to design IBR controls that transform SGs typical second-order response to that of a first-order system. In this way, FSC allows IBRs to provide fast short-term frequency control that fills in the gap produced by the lower system inertia.
However, existing FSC designs are either limited to grid-following (GFL) inverters~\cite{YCVM-2021}, or require IBR buses to be coherent with SG buses~\cite{YJ-AB-PV-EM:21}, a condition that is limited to systems with short electrical distance, i.e., tightly coupled~\cite{MPM-2021}. Moreover, the requirement of shaping SGs as first-order systems leads to high peak power requirements from IBRs. 
Thus, while this approach results in a significant improvement over inertia emulation, its applicability is limited to special settings. 

In this work, we seek to extend this approach to overcome the above-mentioned limitations. Precisely, we consider a setup in which IBRs are electrically distant from SGs (a scenario typically observed in off-shore wind), and provide an FSC design which (i) accounts for the separation of SGs and IBRs, and (ii) trades-off peak power with frequency containment, by matching the frequency dynamics of a low-inertia system to a second-order transfer function. In particular, we define a target second-order response that is motivated by a SG with a smaller effective turbine time constant. Next, we design and analyze a GFM control that realizes this target response under any network coupling. Moreover, we show that stability constraints on the control restrict the range of effective turbine time constants that can be achieved and analyze the robustness of the control towards uncertain network parameters. 

\section{System Modeling}
\label{section:Modeling}
The focus of this work is to investigate shaping the dynamic response of the frequency of a low-inertia power system. For brevity of the presentation, we will consider a two-bus system\footnote{While the network topology impacts system performance and stability \cite{PGD2019}, this aspect is beyond the scope of this manuscript.} (see Figure~\ref{fig:networkforming}) containing a synchronous machine (SM) that models the aggregate frequency dynamics of a conventional multi-machine power system \cite{FP-EM:20} and a grid-forming voltage source inverter (VSI). 

The aggregate response of the SM frequency deviation $\omega_\text{sm}$ to the deviation of the power $p_\text{sm} \in \mathbb{R}$ injected by the SM, is modeled by swing dynamics $G_{\text{sw}}(s)$  with a first-order turbine/governor model $G_{\text{tg}}(s)$, where
\begin{equation*}
G_{\text{sw}}(s)\coloneqq \dfrac{1}{2 H s + \alpha_\ell}, \,\,G_{\text{tg}}(s)\coloneqq -\dfrac{\alpha_{\text{g}}}{\tau s + 1},
\end{equation*}
and results in \cite{FP-EM:20}
\begin{align}\label{eq:machine}
   \omega_{\text{sm}}(s) = -\underbrace{\dfrac{G_{\text{sw}}(s)}{1+G_{\text{sw}}(s) G_{\text{tg}}(s)}}_{\eqqcolon G_{\omega_\text{sm},\,p_\text{sm}}(s)}\, p_{\text{sm}},
\end{align}
where $\alpha_\text{g}\in \mathbb{R}_{>0}$ denotes the aggregate speed governor gain (i.e., inverse frequency droop constant), $\alpha_\ell \in \mathbb{R}_{>0}$ denotes the aggregate frequency sensitivity of load, and $\tau \in \mathbb{R}_{>0}$ is the aggregate turbine time constant. 
\begin{figure}[h]
\centering
\tikzset{every picture/.style={scale=1.0}}%


\usetikzlibrary{calc,trees,positioning,arrows,chains,shapes.geometric,%
    decorations.pathreplacing,decorations.pathmorphing,shapes,%
    matrix,shapes.symbols}

\tikzset{
>=stealth',
  punktchain/.style={
    rectangle, 
    rounded corners, 
    draw=black, very thick,
    text width=10em, 
    minimum height=3em, 
    text centered, 
    on chain},
  line/.style={draw, thick, <-},
  element/.style={
    tape,
    top color=white,
    bottom color=blue!50!black!60!,
    minimum width=8em,
    draw=blue!40!black!90, very thick,
    text width=10em, 
    minimum height=3.5em, 
    text centered, 
    on chain},
  every join/.style={->, thick,shorten >=1pt},
  decoration={brace},
  tuborg/.style={decorate},
  tubnode/.style={midway, right=2pt},
}

\tikzstyle{block} = [draw, fill=white, rounded corners, minimum height=3em, minimum width=6em]
\tikzstyle{sum} = [draw, fill=white, circle, node distance=1cm]
\tikzstyle{input} = [coordinate]
\tikzstyle{output} = [coordinate]
\tikzstyle{pinstyle} = [pin edge={to-,thin,black}]

\begin{circuitikz}[american voltages, american currents]
\draw (0,1.25)  node[ground] {}  to [sV] (0,2.5) node [anchor=east]{\small$p_\text{sm}$};
\draw [-to] (0,2.5)-- (0.25,2.5);
\draw  (0,2.3) to (0, 2.5) to
 [short, -] (2,2.5) to [L, l={\small$b$}] (4,2.5) 
 to  [short, -] (5,2.5);
 \draw
  (5,1.25)  node[ground] {} to [sV] (5,2.5) node [anchor=west]{\small$p_\text{vsi}$};
  \draw [-to] (5,2.5)--(4.75,2.5);
  \draw (5, 2.3)--(5, 2.5);
  \draw  [ultra thick] (1.25,1.75) --(1.25,3)node [anchor=east] {\small$\theta_\text{sm}$} ;
  \draw  [ultra thick] (4,1.75) --(4,3)node [anchor=west] {\small$\theta_\text{vsi}$} ;
  \draw (1.25, 2) to (1.95, 2);
  \draw [-to](1.95, 2) to (1.95, 1.25) node [anchor=north] {\small$p_\ell$};
   \node[left=10pt] at (0,1.375){\small $\omega_\text{sm}$};
 \node[right=10pt] at (5,1.375){\small $\omega_\text{vsi}$};
\end{circuitikz}
\caption{Interconnection of a synchronous machine and grid-forming VSI. 
\label{fig:networkforming}}
\smallskip
\centering
\definecolor{mycolor1}{HTML}{BE9063}%
\definecolor{mycolor2}{RGB}{168,50,45}
\definecolor{mycolor3}{HTML}{525B56}
\definecolor{mycolor4}{HTML}{A4978E}%
\definecolor{mycolor5}{rgb}{0.59400,0.18400,0.35600}%
\definecolor{mycolor6}{HTML}{253F5B}%
\definecolor{mycolor7}{HTML}{818A6F}
\definecolor{mycolor8}{HTML}{D59B2D}
\definecolor{mycolor9}{RGB}{31,140,24}
\definecolor{mycolor10}{RGB}{136,176,197}

\newcommand{\tb}[1]{\textcolor{blue}{#1}}
\newcommand{\tg}[1]{\textcolor{ForestGreen}{#1}}
\newcommand{\tol}[1]{\textcolor{OliveGreen}{#1}}
\newcommand{\tor}[1]{\textcolor{Orange}{#1}}
\newcommand{\tnb}[1]{\textcolor{NavyBlue}{#1}}

\usetikzlibrary{calc,trees,positioning,arrows,chains,shapes.geometric,%
    decorations.pathreplacing,decorations.pathmorphing,shapes,%
    matrix,shapes.symbols}

\tikzset{
>=stealth',
  punktchain/.style={
    rectangle, 
    rounded corners, 
    draw=black, very thick,
    text width=10em, 
    minimum height=3em, 
    text centered, 
    on chain},
  line/.style={draw, thick, <-},
  element/.style={
    tape,
    top color=white,
    bottom color=blue!50!black!60!,
    minimum width=8em,
    draw=blue!40!black!90, very thick,
    text width=10em, 
    minimum height=3.5em, 
    text centered, 
    on chain},
  every join/.style={->, thick,shorten >=1pt},
  decoration={brace},
  tuborg/.style={decorate},
  tubnode/.style={midway, right=2pt},
}

\tikzstyle{block} = [draw, fill=white, rounded corners, minimum height=3em, minimum width=6em]
\tikzstyle{sum} = [draw, fill=white, circle, minimum size=0.6cm, node distance=1cm]
\tikzstyle{input} = [coordinate]
\tikzstyle{output} = [coordinate]
\tikzstyle{pinstyle} = [pin edge={to-,thin,black}]

\begin{tikzpicture}[auto, node distance=2cm,>=latex']

    \node [input, name=input] {};
    \node [sum, right=0.75cm of input] (sum) {};
    \node [block, right=1cm of sum, fill=mycolor4] (controller) {\small $G_{\omega_\text{sm},\, p_\text{sm}}(s)$};
    \node [output, right=1.05cm of controller] (temp) {};
    \node [output, right=0.75cm of temp] (output) {};
     \node [block, below=0.5cm of controller, fill=mycolor10] (turbine) {\small$G_{\omega_\text{vsi}, \, p_\text{vsi}}(s)$};
       \node[left = 1.3cm of turbine] (diff){};
       \node[sum, right = 0.75cm of turbine] (diff2){};
      \node [block, below= 0.5cm of turbine, fill=mycolor10] (IBR) {\small$\dfrac{1}{s}\,b$};
        \node[ below= 0.17cm of turbine](disptemp){};
   \draw [-] (controller) -- node {} (temp);
    \draw [->] (temp) -- node {\small$\omega_\text{sm}$} (output);
    \draw [->] (sum) -- node {} (controller);
     \draw [<-] (diff2) -- node {\small$\omega_\text{vsi}$} (turbine);
     \draw [<-] (diff2) -- node {} (temp);
     \draw [->] (diff2) |- node {} (IBR);
     \draw [->] (IBR) -| node[pos=0.95] {\small$p_\text{vsi}$} (sum);
    \draw [->] (input) -- node {\small{$p_\ell$}} (sum);
    \draw [<-] (turbine) -- (diff);

    \draw [-, dashed, very thick, draw=mycolor1] (0.4,-0.8) -- (0.4, -4);   
    \draw [-,  dashed, very thick, draw=mycolor1] (0.4,-4) -- (6, -4); 
    \draw [-,  dashed, very thick, draw=mycolor1] (6,-4) -- (6, -0.8); 
    \draw [-,  dashed, very thick, draw=mycolor1] (6, -0.8) -- (0.4, -0.8);

    \draw (sum.north east) -- (sum.south west)
    (sum.north west) -- (sum.south east);
    
    \draw (diff2.north east) -- (diff2.south west)
    (diff2.north west) -- (diff2.south east);
    
    \node[left=-0.25pt] at (sum.center){\tiny$-$};
    \node[below=0.25pt] at (sum.center){\tiny $-$};
    \node[left=-0.25pt] at (diff2.center){\tiny$-$};
    \node[above=0.25pt] at (diff2.center){\tiny $+$};
    \node[text=mycolor1] at (1.35,-3.7) {\small$-G_{p_\text{vsi}, \, \omega_\text{sm}}(s)$};



\end{tikzpicture}
\caption{The synchronous machine--IBR viewed as a feedback system.
\label{fig:IBR-machine}}
\end{figure}

Next, consider a VSI coupled to the aggregate frequency dynamics \eqref{eq:machine} through a lossless transmission line with susceptance $b \in \mathbb{R}_{>0}$ as depicted in Figure~\ref{fig:networkforming}. We will refer to the grid as weakly-coupled if $b$ is small (e.g., for grids with low short-circuit ratio) and tightly-coupled if $b$ is large. Notably, this relaxes the assumption in \cite{YJ-AB-PV-EM:21} that $b\to \infty$. 
Next, let $\theta_\text{vsi}(s)=\tfrac{1}{s} \omega_\text{vsi}(s)$ and $\theta_\text{sm}(s)=\tfrac{1}{s}\omega_\text{sm}(s)$ denote the voltage phase angles at the IBR and SM bus. Then, using the DC power flow approximation at $1~\mathrm{p.u.}$ voltage magnitude and zero angle difference \cite{FP-EM:20}, the IBR power injection $p_\text{vsi}(s)$ is equal to the power flowing across the line, i.e.,
\begin{align}\label{eq:networktf}
    p_\text{vsi}(s)= b\,\big(\theta_\text{vsi}(s)-\theta_{\text{sm}}(s)\big)=\dfrac{1}{s}\,b\, \big(\omega_\text{vsi}(s)-\omega_{\text{sm}}(s)\big).
\end{align}
While the quasi-steady-state model \eqref{eq:networktf} cannot detect instability due to network circuit dynamics, its use for the problem at hand is  justified by large transformer impedances in transmission systems \cite[Fig.~4]{G2022} and the control developed in Sec.~\ref{section:IBR_implementation}.

The combination of the synchronous machine and the IBR can be interpreted as the VSI in feedback with the synchronous machine as shown in Figure~\ref{fig:IBR-machine}. 
To this end, let
\begin{align}\label{eq:vsitf}
\omega_{\text{vsi}}(s)=-G_{\omega_\text{vsi},\, p_\text{vsi}}(s)\, p_\text{vsi}(s),
\end{align}
represent the dynamics of the grid-forming VSI. Combing \eqref{eq:networktf} and \eqref{eq:vsitf}, the relation between the generator frequency $\omega_{\text{sm}}$ and the inverter power $p_\text{vsi}(s)$ is
\begin{align}
\label{eq:IBRinjectionP}
    p_\text{vsi}(s)=\underbrace{-\dfrac{1}{\tfrac{1}{b}s+G_{\omega_\text{vsi},\, p_\text{vsi}}(s)}}_{\textstyle\coloneqq G_{p_\text{vsi}, \, \omega_\text{sm}}(s)}\,\omega_{\text{sm}}(s).
\end{align}
Using the load perturbation $p_\ell$, $p_{\text{sm}}=p_\ell - p_\text{vsi}$, and \eqref{eq:IBRinjectionP} to close the loop between \eqref{eq:vsitf} and \eqref{eq:machine}, results in the \emph{closed-loop} transfer function
\begin{align}
    \omega_{\text{sm}}(s)=\underbrace{-\dfrac{(\tfrac{1}{b}s+G_{\omega_\text{vsi},\, p_\text{vsi}}(s))\,G_{\omega_\text{sm},\, p_\text{sm}}(s)}{\tfrac{1}{b}s+G_{\omega_\text{sm},\, p_\text{sm}}(s)+G_{\omega_\text{vsi},\, p_\text{vsi}}(s)}}_{\textstyle \coloneqq G^{\text{cl}}_{\omega_\text{sm}, \, p_\ell}(s)}\, p_\ell(s),
    \label{eq:overallresp}
\end{align}
from the load perturbation $p_\ell$ to the aggregate frequency $\omega_\text{sm}$. 

\section{Problem Formulation}
\label{Sec:Prob}
The goal of this work is to shape the response $G^{\text{cl}}_{\omega_\text{sm}, \, p_\ell}(s)$ through the IBR control $G_{\omega_\text{vsi},\, p_\text{vsi}}(s)$ to improve the system frequency response. To that end, in this section, we will first discuss the desired response $G^{\text{cl}\star}_{\omega_\text{sm}, \, p_\ell}(s)$. Subsequently, in Section~\ref{section:IBR_implementation}, we will compute $G_{\omega_\text{vsi},\, p_\text{vsi}}(s)$ such that $G^{\text{cl}}_{\omega_\text{sm}, \, p_\ell}(s)$ matches the desired response $G^{\text{cl}\star}_{\omega_\text{sm}, \, p_\ell}(s)$. In \cite{YJ-AB-PV-EM:21}, the target transfer function of the overall system $G^{\text{cl}\star}_{\omega_\text{sm}, \, p_\ell}(s)$ is a \emph{first-order} system. As first-order responses do not exhibit overshoots, the post-disturbance frequency nadir and the steady-state settling frequency are identical. This results in the corresponding IBR response to exhibit a high peak power injection $p_\text{vsi}(s)$ in order to reduce the frequency nadir. Power inverters with high rating (i.e., that can provide large peak power injections) are expensive and frequent large power injections degrade the lifetime, e.g., battery energy storage systems. Thus, it is critical to minimize the IBR peak power, without potentially degrading the frequency response of the overall system too much. In particular, we aim to fully leverage the fast and flexible actuation capabilities of IBRs to improve the frequency nadir, while lowering the IBR peak power injections relative to the approach in \cite{YJ-AB-PV-EM:21}. Next, we illustrate that the power injections from the IBRs can be improved by choosing \emph{second-order} target transfer functions $G^{\text{cl}\star}_{\omega_\text{sm}, \, p_\ell}(s)$. 

\subsection{Generator matching transfer function}
We consider the standalone synchronous machine transfer function $G_{\omega_\text{sm},\,p_\text{sm}}(s)$ defined in \eqref{eq:machine}. Next, we define the set 
\begin{align}
\mathcal{U}\coloneqq \left\{\rho \in \mathbb{R} \,\bigg\vert\, 0 \leq \rho < \tau \right\}
\end{align}  
of effective turbine time constants. As we desire to design a second-order transfer function for our target response $G^{\text{cl}\star}_{\omega_\text{sm}, \, p_\ell}(s)$, we consider a candidate transfer function in the form of \eqref{eq:machine}  with an effective turbine time constant $\rho \in \mathcal{U}$, i.e., 
\begin{align}
    G^{\text{cl}\star}_{\omega_\text{sm}, \, p_\ell}(s) \coloneqq -\frac{s\,\rho+1}{2H\rho \,s^2+(\alpha_\ell\rho+2H)\,s+(\alpha_\ell+\alpha_g)}.
    \label{eq:targetresptime}
\end{align}
Such a choice of the target response translates to the IBR speeding up the response of the turbine/governor system of the SG, i.e., an overall system response that is much faster than the original (aggregate) power system \eqref{eq:machine}. This choice is justified by the fact that SGs provide significant (but slow) frequency control, whereas IBRs can provide a fast frequency response but often cannot sustain this response over long periods of time due to limited flexibility and energy storage.

\begin{remark}\longthmtitle{Target second-order response}
{\rm In our analysis, we consider the target $G^{\text{cl}\star}_{\omega_\text{sm}, \, p_\ell}(s)$ to be a modification of the original system response $G_{\omega_\text{sm},\,p_\text{sm}}(s)$. Such a choice is reasonable as the IBR devices cannot provide significant levels of inertia without significantly oversizing the IBR \cite{TGA+20}. Moreover, we aim to control the post-disturbance steady-state power injection of the IBR to zero. These requirements effectively prevent us from modifying the inertia constant $H$, load damping $\alpha_\ell$, and governor gain $\alpha_g$. Considering more generic target transfer functions (parameterized by poles, damping ratio, and zeros)  that allow specifying objectives beyond minimizing nadir and peak power injection are seen as an interesting area for future work.
}
\oprocend
\end{remark}

In general, we can design controllers for IBRs to realize a wide range of power injections. In order to  obtain an overall system response $G^{\text{cl}}_{\omega_\text{sm}, \, p_\ell}(s)=G^{\text{cl}\star}_{\omega_\text{sm}, \, p_\ell}(s)$, the power injection from the IBR, $p_\text{vsi}(s)$ has to satisfy
\begin{align}
\label{eq:IBRtarget}
 p_\text{vsi}(s)\overset!= -&\dfrac{\alpha_g\,s\, (\tau-\rho)}{(s\,\tau +1)(s\,\rho+1)}\,\omega_{\text{sm}}(s),
 \end{align}
 where $\omega_{\text{sm}}(s)$ is the frequency of the synchronous machine \eqref{eq:machine}. In terms of the load disturbance $p_\ell(s)$, the above relation in conjunction with \eqref{eq:overallresp} translates to
 \begin{align}
 p_\text{vsi}(s)=- &\dfrac{\alpha_g\,s\, (\tau-\rho)}{(s\,\tau +1)(s\,\rho+1)}\,G^{\text{cl}\star}_{\omega_\text{sm}, \, p_\ell}(s)\,p_\ell(s),
 \label{eq:powerIBRload}
 \end{align}
 where the target response $G^{\text{cl}\star}_{\omega_\text{sm}, \, p_\ell}(s)$ is given by \eqref{eq:targetresptime}. We note from \eqref{eq:powerIBRload} that the transfer function is third-order due to stable pole-zero cancellation. 

\begin{remark}\longthmtitle{IBR power injection}
{\rm While we limit ourselves to grid-forming inverters in this paper, the results in this section are  agnostic to the IBR implementation. As long as an IBR injects the power specified by \eqref{eq:powerIBRload}, the overall system response will match the target $G^{\text{cl}\star}_{\omega_\text{sm}, \, p_\ell}(s)$.}
\oprocend
\end{remark}
 
 \subsection{Minimizing Peak IBR power injections}
 We recall from the discussion above that the main motivation for exploring higher-order target transfer functions is to minimize the peak power injected by the IBR devices while limiting the aggregate system frequency excursion in response to load perturbations. Let $p_\text{vsi}(t)$ be the power injected by the IBR in response to a step perturbation in the load $p_\ell$, and $|\cdot|_\infty=\sup_t|\cdot(t)|$ denote the absolute peak value, then the problem of minimizing peak IBR power translates to
 \begin{align}
 \label{eq:originalprob}
     \min_{\rho} \quad & |p_\text{vsi}|_\infty \qquad
     \text{s.t.}\quad |\omega_{\text{sm}}|_\infty \leq \bar{\omega}_{\text{sm}}.
 \end{align}
 Including IBR current limits and power reserve capacity constraints is an interesting area for future research.
 
 While \cite{FP-EM:20} provides a closed-form relation for the frequency nadir of an under-damped SG of \eqref{eq:machine}, exact solutions to the problem \eqref{eq:originalprob} are generally intractable both analytically and computationally. Thus, we resort to gridding the parameter space. Figure~\ref{fig:freqresptime} and Figure~\ref{fig:powerinjectime} depict the  response of the SG frequency and IBR power injection to a $1~\mathrm{p.u.}$  load step for different overall effective time constants $\rho$. We note that as the effective time constant $\rho$ is reduced, the corresponding frequency nadir decreases due to higher IBR power injection.

\begin{figure}[h!]
\centering
\tikzset{every picture/.style={scale=1.0}}%
\input{freq_resp.tex}
\caption{The frequency response for a $1~\mathrm{p.u.}$ load step for different $\rho$.}
\label{fig:freqresptime}
\medskip
\centering
\tikzset{every picture/.style={scale=1.0}}%
\input{IBR_power.tex}
\caption{The IBR power injection for a $1~\mathrm{p.u.}$ load step for different $\rho$.}
\label{fig:powerinjectime}
\end{figure}

Another observation for the specific system parameters used\footnote{We use the single machine-single IBR system from \cite{YJ-AB-PV-EM:21} for our simulations.}, is that the objective is monotone in the decision variable $\rho$ and thus the solution to the minimization problem \eqref{eq:originalprob} can also be determined graphically. To this end, we compute and plot specifically, the peak power injection from IBRs for step load perturbations in Figure~\ref{fig:time-nadir} as a function of $\rho$. We also plot the frequency nadir\footnote{We consider the absolute value of the frequency nadir here. As the load perturbation is a step, the frequency nadir  will be negative.} of the overall system (with the IBR) as a function of the effective time constant $\rho$ based on the data from \cite{YJ-AB-PV-EM:21}. As the effective time constant increases, the frequency nadir increases monotonically attaining a maximum for $\rho=\tau$, i.e., the original system without any IBR. Finally, in Figure~\ref{fig:paretop}, we plot the Pareto front for frequency nadir and the peak power input $p_\text{vsi}$ from the IBR. We also indicate points on the Pareto front for varying $\rho$. The limit (i.e., $\rho\to 0$) recovers the first-order response as in \cite{YJ-AB-PV-EM:21}. Note that as both these quantities are monotonic in $\rho$, the Pareto front would allow a system operator to  choose an acceptable trade-off between the two.
\begin{figure}[htbp]
\centering
\tikzset{every picture/.style={scale=1.0}}%
\definecolor{mycolor1}{HTML}{BE9063}%
\definecolor{mycolor2}{RGB}{168,50,45}
\definecolor{mycolor3}{HTML}{525B56}
\definecolor{mycolor4}{HTML}{A4978E}%
\definecolor{mycolor5}{rgb}{0.59400,0.18400,0.35600}%
\definecolor{mycolor6}{HTML}{253F5B}%
\definecolor{mycolor7}{HTML}{818A6F}
\definecolor{mycolor8}{HTML}{D59B2D}
\definecolor{mycolor9}{RGB}{31,140,24}
\definecolor{mycolor10}{RGB}{136,176,197}

\begin{tikzpicture}
\pgfplotsset{
    scale only axis,
    xmin=0.1, xmax=1
}
\begin{axis}[
width=2.6in,
height=1.2in,
at={(1.739in,0.849in)},
  axis y line*=left,
  ymin=200, ymax=450,
 xticklabel style = {font=\footnotesize,yshift=0ex},
xlabel style={font=\color{black}},
xlabel={\small Effective time constant $\rho$ (s)},
  ymajorgrids,
  xmajorgrids,
yticklabel style = {font=\footnotesize,xshift=0ex},
ylabel style={font=\color{black}},
ylabel={\small Frequency nadir (mHz)},
 xtick={0.1, 0.25, 0.4, 0.55, 0.7, 0.85, 1.0},
 legend style={legend cell align=left, align=left, draw=black, font=\small, draw=none, legend columns=1, at={(0.6,0.85)}}
]
\addplot[smooth,color=mycolor1, line width=1.75pt]
  coordinates{
   (0.1,	211.215244610773)
(0.15,	214.569967689315)
(0.2,	237.099422510742)
(0.25,	255.483677460081)
(0.3,	271.850407845722)
(0.35,	286.837130248348)
(0.4,	300.758677900104)
(0.45,	313.811576797325)
(0.5,	326.13346596622)
(0.55,	337.827112446415)
(0.6,	348.972412047316)
(0.65,	359.633261333538)
(0.7,	369.86186133479)
(0.75,	379.701584037251)
(0.8,	389.18896868097)
(0.85,	398.355159940232)
(0.9,	407.226972173342)
(0.95,	415.827694416735)
(1,	424.177710580167)
}; \addlegendentry{nadir}
\end{axis}

\begin{axis}[
width=2.6in,
height=1.1in,
at={(1.739in,0.849in)},
  axis y line*=right,
  axis x line=none,
  ymin=0, ymax=0.6,
yticklabel style = {font=\footnotesize,xshift=0ex},
ylabel style={font=\color{black}},
ylabel={\small Peak IBR power (p.u.)},
 legend style={legend cell align=left, align=left, draw=black, font=\small, draw=none, legend columns=1, at={(0.7,1.08)}}
]

\addplot[smooth,color=mycolor2, line width=1.75pt]
   coordinates{
(0.1,	0.580209654966924)
(0.15,	0.567773220643045)
(0.2,	0.545341275691656)
(0.25,	0.517239860610336)
(0.3,	0.485993769840851)
(0.35,	0.452721742230255)
(0.4,	0.418307175705351)
(0.45,	0.383119595376276)
(0.5,	0.34744769893605)
(0.55,	0.311508776662964)
(0.6,	0.27558446625031)
(0.65,	0.24016316090152)
(0.7,	0.204772033036198)
(0.75,	0.169641481637973)
(0.8,	0.135018331344184)
(0.85,	0.100621046725119)
(0.9,	0.0667065503472208)
(0.95,	0.0331451324789037)
(1,	7.61307890438353e-17)
}; \addlegendentry{peak power}
\end{axis}

\end{tikzpicture}
\caption{Frequency nadir and peak power injection for a $1~\mathrm{p.u.}$ load step as a function of $\rho$.}\label{fig:time-nadir}
\medskip
  \centering
\tikzset{every picture/.style={scale=1.0}}%
%
%
\definecolor{mycolor1}{HTML}{BE9063}%
\definecolor{mycolor2}{RGB}{168,50,45}
\definecolor{mycolor3}{HTML}{525B56}
\definecolor{mycolor4}{HTML}{A4978E}%
\definecolor{mycolor5}{rgb}{0.59400,0.18400,0.35600}%
\definecolor{mycolor6}{HTML}{253F5B}%
\definecolor{mycolor7}{HTML}{818A6F}
\definecolor{mycolor8}{HTML}{D59B2D}
\definecolor{mycolor9}{RGB}{31,140,24}
\definecolor{mycolor10}{RGB}{136,176,197}

\begin{tikzpicture}
\begin{axis}[%
width=2.9in,
height=1.2in,
at={(1.739in,0.849in)},
scale only axis,
xmin=200,
xmax=450,
ymin=0,
ymax=0.6,
ytick={0, 0.15, 0.3, 0.45, 0.6, 0.75},
ymajorgrids,
xmajorgrids,
yticklabel style = {font=\footnotesize,xshift=0ex},
xticklabel style = {font=\footnotesize,yshift=0ex},
ylabel style={font=\color{black}},
ylabel={\small Peak IBR power (p.u.)},
xlabel style={font=\color{black}},
xlabel={\small Frequency nadir (mHz)},
axis background/.style={fill=white},
legend style={legend cell align=left, align=left, draw=black, font=\small, draw=none, legend columns=-1, at={(1.02,1.22)}}
]
\addplot [color=mycolor1, line width=1.75pt]
  table[row sep=crcr]{%
211.215244610773	0.580209654966924\\
214.569967689315	0.567773220643045\\
237.099422510742	0.545341275691656\\
255.483677460081	0.517239860610336\\
271.850407845722	0.485993769840851\\
286.837130248348	0.452721742230255\\
300.758677900104	0.418307175705351\\
313.811576797325	0.383119595376276\\
326.133465966226	0.347447698936057\\
337.827112446415	0.311508776662964\\
348.972412047316	0.27558446625031\\
359.633261333538	0.24016316090152\\
369.861861334796	0.204772033036198\\
379.701584037251	0.169641481637973\\
389.18896868097	0.135018331344184\\
398.355159940232	0.100621046725119\\
407.226972173342	0.0667065503472208\\
415.827694416735	0.0331451324789037\\
424.177710580167	7.61307890438353e-17\\
};

\addplot[only marks, mark=*, mark options={}, mark size=2pt, draw=mycolor1, fill=mycolor2] table[row sep=crcr]{%
x	y\\
424.177710580167		7.61307890438353e-17\\
};

\addplot[only marks, mark=*, mark options={}, mark size=2pt, draw=mycolor1, fill=mycolor2] table[row sep=crcr]{%
x	y\\
389.18896868097	0.135018331344184\\
};

\addplot[only marks, mark=*, mark options={}, mark size=2pt, draw=mycolor1, fill=mycolor2] table[row sep=crcr]{%
x	y\\
348.972412047316	0.27558446625031\\
};

\addplot[only marks, mark=*, mark options={}, mark size=2pt, draw=mycolor1, fill=mycolor2] table[row sep=crcr]{%
x	y\\
300.758677900104	0.418307175705351\\
};

\addplot[only marks, mark=*, mark options={}, mark size=2pt, draw=mycolor1, fill=mycolor2] table[row sep=crcr]{%
x	y\\
237.099422510742	0.545341275691656\\
};

\node at ($(190,35)$) {\small no IBR};
\node at ($(163,130)$) {\small $\rho=0.8$};
\node at ($(127,250)$) {\small $\rho=0.6$};
\node at ($(80,390)$) {\small $\rho=0.4$};
\node at ($(25,490)$) {\small $\rho=0.2$};

\end{axis}
\end{tikzpicture}%
 \caption{Pareto front for peak IBR power injection v/s the frequency nadir for a $1~\mathrm{p.u.}$ load step. The limit $\rho\to 0$ recovers a first-order response.}
\label{fig:paretop}
\end{figure}

Thus, for the target system response in \eqref{eq:targetresptime}, we obtain a set of pairs of peak IBR power injection and frequency nadir points. Each pair maps to a unique $\rho$, which has the interpretation of time constant analogous to the turbine time constant of the original system \eqref{eq:machine}. 

\section{IBR implementation of Controllers}
\label{section:IBR_implementation}

In this section, we delve into the implementation of the frequency controls discussed in the previous section using IBRs. In particular, we investigate grid-forming devices, which we believe will form the backbone of future grids. 

Next, consider the single machine, single grid-forming IBR system as in Figure~\ref{fig:networkforming}. We emphasize that, compared to the case $b\to \infty$ considered in \cite{YJ-AB-PV-EM:21}, both  weak grid-coupling (i.e., small $b$) and tight grid-coupling (i.e., large but finite $b$) can have a non-trivial impact on frequency shaping control.

Let $\theta_{\text{sm}}$, $p_\text{sm}$ (resp. $\theta_\text{vsi}$, $p_\text{vsi}$) denote the angle, power injection of the synchronous machine (resp. grid-forming inverter) and $p_\ell$ denote the load served by the system. The line susceptance (which may be time-varying, e.g., depending on the state of tap changing transformers, etc.,) between the IBR and the machine is modeled by the susceptance $b$. For this system, let  the synchronous machine dynamics be as in \eqref{eq:machine}, and the inverter dynamics as in \eqref{eq:networktf}. We wish to design the control transfer function $G_{\omega_\text{vsi},\, p_\text{vsi}}(s)$ in order to realize the target function $G^{\text{cl}\star}_{\omega_\text{sm}, \, p_\ell}(s)$ for the overall system. We recall from \eqref{eq:IBRtarget} the IBR power injection required to achieve the target response, i.e., $G^{\text{cl}}_{\omega_\text{sm}, \, p_\ell}(s)=G^{\text{cl}\star}_{\omega_\text{sm}, \, p_\ell}(s)$. On comparing with \eqref{eq:IBRinjectionP}, we obtain
\begin{align}
{G^{-1}_{p_\text{vsi}, \, \omega_\text{sm}}(s)}=\dfrac{1}{b}\,s+G_{\omega_\text{vsi},\, p_\text{vsi}}(s)\overset!=\dfrac{(s\,\tau +1)(s\,\rho+1)}{\alpha_g\,s\, (\tau-\rho)}.
\label{eq:PID_matching}
\end{align}

We note that the right hand side terms of \eqref{eq:PID_matching} can be realized through a Proportional-Integral-Derivative (PID) type controller $G_{\omega_\text{vsi},\, p_\text{vsi}}(s)$, i.e.,
\begin{align}
G_{\omega_\text{vsi},\, p_\text{vsi}}(s)=k_p + \dfrac{k_i}{s}+ k_d\, s.
\end{align}

We note that the derivative gain $k_d$ can be used to compensate network circuit dynamics \cite[Sec.~V]{G2022} and thereby justify the use of the quasi-steady-state network model \eqref{eq:networktf} for stability analysis. Next, we consider the PID control gains
\begin{align}
\label{eq:pid_gains}
\!\! k_d \!\coloneqq\! \dfrac{\tau\, \rho}{\alpha_g (\tau-\rho)}\!-\!\dfrac{1}{\hat{b}},\, k_p \!\coloneqq\! \dfrac{\tau+ \rho}{\alpha_g (\tau-\rho)},\,
k_i \!\coloneqq\! \dfrac{1}{\alpha_g (\tau-\rho)},\!\!
\end{align}
where $\hat{b}$ is an estimate for $b$. If $\hat{b}=b$, then \eqref{eq:PID_matching} holds and closed-loop transfer $G^{\text{cl}}_{\omega_\text{sm}, \, p_\ell}$ function equals \eqref{eq:targetresptime}. We will first analyze the case $\hat{b}=b$. To this end, the set of effective time constant for which $k_d\geq 0$ is denoted by
\begin{align}
\mathcal{N}\coloneqq \left\{\rho \in \mathbb{R} \,\bigg\vert\, \dfrac{\alpha_g\, \tau}{\hat{b}\, \tau+\alpha_g} \leq \rho < \tau \right\} \subset \mathcal{U}.
\end{align}

\begin{proposition}
\label{propo:stablesys}
Consider the control gains $k_d$, $k_p$, and $k_i$ defined in \eqref{eq:pid_gains}, and $\hat{b}=b$, then
\begin{enumerate}[label=(\roman*)]
    \item the transfer function $G_{\omega_\text{sm},\, p_\text{vsi}}(s)$ is stable,
    \item the closed-loop transfer function $G^{\text{cl}}_{\omega_\text{sm}, \, p_\ell}(s)$ is passive,
    \item and $-G_{\omega_\text{vsi},\, p_\text{vsi}}(s)$ is minimum phase iff $\rho \in \mathcal{N}$.
\end{enumerate}
\end{proposition}
\begin{proof}
For the PID control gains in \eqref{eq:pid_gains}, the resulting transfer function $G_{\omega_\text{sm},\, p_\text{vsi}}(s)$ is stable from the Routh-Hurwitz criterion. Next, we note that the effective turbine dynamics
\begin{align}
    G_{\omega_\text{sm}, \, p_\text{vsi}}(s)-\dfrac{\alpha_g}{s\,\tau + 1} = -\dfrac{\alpha_g}{s\,\rho + 1},
\end{align}
is passive. Moreover, the overall system transfer function $G^{\text{cl}}_{\omega_\text{sm}, \, p_\ell}(s)$ is a negative feedback interconnection of the strictly passive machine rotor dynamics and the passive effective turbine dynamics, thus stable. 
 The last item directly follows from $k_d\geq0$.
\end{proof}
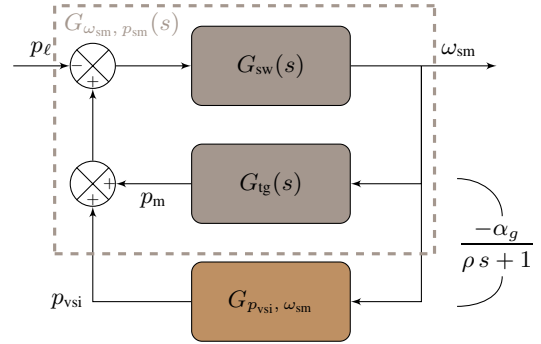
\begin{figure}[htbp]
\centering
\tikzset{every picture/.style={scale=1.0}}%
\definecolor{mycolor1}{HTML}{BE9063}%
\definecolor{mycolor2}{RGB}{168,50,45}
\definecolor{mycolor3}{HTML}{525B56}
\definecolor{mycolor4}{HTML}{A4978E}%
\definecolor{mycolor5}{rgb}{0.59400,0.18400,0.35600}%
\definecolor{mycolor6}{HTML}{253F5B}%
\definecolor{mycolor7}{HTML}{818A6F}
\definecolor{mycolor8}{HTML}{D59B2D}
\definecolor{mycolor9}{RGB}{31,140,24}
\definecolor{mycolor10}{RGB}{136,176,197}

\newcommand{\tb}[1]{\textcolor{blue}{#1}}
\newcommand{\tm}[1]{\textcolor{magenta}{#1}}
\newcommand{\tg}[1]{\textcolor{ForestGreen}{#1}}
\newcommand{\tol}[1]{\textcolor{OliveGreen}{#1}}
\newcommand{\tor}[1]{\textcolor{Orange}{#1}}
\newcommand{\tnb}[1]{\textcolor{NavyBlue}{#1}}

\usetikzlibrary{calc,trees,positioning,arrows,chains,shapes.geometric,%
    decorations.pathreplacing,decorations.pathmorphing,shapes,%
    matrix,shapes.symbols}

\tikzset{
>=stealth',
  punktchain/.style={
    rectangle, 
    rounded corners, 
    draw=black, very thick,
    text width=10em, 
    minimum height=3em, 
    text centered, 
    on chain},
  line/.style={draw, thick, <-},
  element/.style={
    tape,
    top color=white,
    bottom color=blue!50!black!60!,
    minimum width=8em,
    draw=blue!40!black!90, very thick,
    text width=10em, 
    minimum height=3.5em, 
    text centered, 
    on chain},
  every join/.style={->, thick,shorten >=1pt},
  decoration={brace},
  tuborg/.style={decorate},
  tubnode/.style={midway, right=2pt},
}

\tikzstyle{block} = [draw, fill=white, rounded corners, minimum height=3em, minimum width=6em]
\tikzstyle{sum} = [draw, fill=white, circle, minimum size=0.6cm, node distance=1cm]
\tikzstyle{input} = [coordinate]
\tikzstyle{output} = [coordinate]
\tikzstyle{pinstyle} = [pin edge={to-,thin,black}]

\begin{tikzpicture}[auto, node distance=2cm,>=latex']

    \node [input, name=input] {};
    \node [sum, right=0.75cm of input] (sum) {};
    \node [block, right=1cm of sum, fill=mycolor4] (controller) {\small $G_\text{sw}(s)$};
    \node [output, right of=controller] (temp) {};
        \node [output, right=1cm of temp] (output) {};
     \node [block, below=0.5cm of controller, fill=mycolor4] (turbine) {\small$G_\text{tg}(s)$};
       \node[sum, left = 1cm of turbine] (diff){};
      \node [block, below= 0.5cm of turbine, fill=mycolor1] (IBR) {\small$G_{p_\text{vsi},\, \omega_\text{sm}}$};
        \node[ below= 0.17cm of turbine](disptemp){};
      \node[ right= 2.25cm of disptemp](disp){\small${\dfrac{\small-\alpha_g}{\small\rho\,s+1}}$};
   \draw [-] (controller) -- node {} (temp);
    \draw [->] (temp) -- node {\small$\omega_\text{sm}$} (output);
      \draw [->] (sum) -- node {} (controller);
         \draw [->] (temp) |- node {} (turbine);
            \draw [->] (temp) |- node {} (IBR);
   \draw [->] (turbine) -- node {\small$p_\text{m}$} (diff);
      \draw [->] (IBR) -| node {\small$p_\text{vsi}$} (diff);
            \draw [->] (input) -- node {\small{$p_\ell$}} (sum);
              \draw [->] (diff) -- node {} (sum);
              
              \draw[-] (5.9, -1.5) to[out=0,in=90] (6.5, -2);
               \draw[-] (5.9, -3.2) to[out=0,in=-90] (6.5, -2.8);
              
         \draw [-, dashed, very thick, draw=mycolor4] (0.55,0.8) -- (0.55, -2.5);   
          \draw [-,  dashed, very thick, draw=mycolor4] (0.55,-2.5) -- (5.6, -2.5); 
           \draw [-,  dashed, very thick, draw=mycolor4] (5.6,-2.5) -- (5.6, 0.8); 
            \draw [-,  dashed, very thick, draw=mycolor4] (5.6, 0.8) -- (0.55, 0.8); 
                \node[text=mycolor4] at (1.45,0.55) {\small$G_{ \omega_\text{sm},\, p_\text{sm}}(s)$};
            
            \draw (sum.north east) -- (sum.south west)
    (sum.north west) -- (sum.south east);
    
    \draw (diff.north east) -- (diff.south west)
    (diff.north west) -- (diff.south east);
    
    \node[left=-0.25pt] at (sum.center){\tiny$-$};
\node[below=0.25pt] at (sum.center){\tiny $+$};

\node[right=0.25pt] at (diff.center){\tiny $+$};
\node[below=0.25pt] at (diff.center){\tiny $+$};


\end{tikzpicture}
\caption{The net effect of the grid-forming IBR on the overall system response.}
\label{fig:inv_interconnection}
\end{figure}
Furthermore, we note that the lower bound for the differential gain can be reduced further with the exception that such a controller is non-minimum phase. This choice while improving the allowable range of $\rho$, i.e., $\rho \in \mathcal{U}$ (note that $\mathcal {N} \subset \mathcal{U}$) may however, result in an unstable closed loop if $b$ is not known exactly. We shall investigate this scenario in the next section. Finally, we remark that the methods developed in this section can be applied to more general topologies.

\begin{remark}\longthmtitle{IBRs in star topology}
\rm The closed-loop transfer function \eqref{eq:overallresp} can be generalized for multiple IBRs connected to the synchronous generator in a star topology as
\begin{align*}
    \omega_{\text{sm}}(s)=-\dfrac{G_{\omega_\text{sm},\, p_\text{sm}}(s)}{1+\sum_i \dfrac{G_{\omega_\text{sm},p_\text{sm}}(s)}{\tfrac{s}{b_i}+G_{\omega_\text{vsi,i}, \, p_\text{vsi,i}}(s)}}\, p_\ell(s),
\end{align*}
where $b_i$ denotes the line susceptance of the $i$-th IBR. For this configuration, the overall response from the IBRs can be distributed, such that the sum of the contributions equals \eqref{eq:PID_matching}, 
\begin{align*}
\sum_i\bigg(\dfrac{s}{b_i}+G_{\omega_\text{vsi,i},\, p_\text{vsi,i}}(s)\!\bigg)^{-1}\overset!=\dfrac{\alpha_g\,s\, (\tau-\rho)}{(s\,\tau +1)(s\,\rho+1)},
\end{align*}
and this is a compelling extension for future research. 
\oprocend
\end{remark}

\section{Effect of imperfect network parameters}
In the preceding discussions, we have assumed that the network susceptance is known with a fair degree of certainty. However, in various practical settings the network susceptance may generally be unknown, may be time-varying, and potentially can only be approximately estimated. For such scenarios, analyzing the robustness and performance trade-offs merits a detailed discussion. Moreover, from \eqref{eq:pid_gains} we note that the parameters $\alpha_\ell$ and $H$ do not impact the control. While $\alpha_g$ impacts the control gains, it is explicitly known from, e.g., grid codes or markets in many systems. Thus, we only focus on the impact of the network susceptance $b$.

\subsection{Stability and sensitivity analysis}
We will now consider controller gains that are functions of the estimate $\hat{b}\neq b$. The analysis in the previous section crucially hinges on a pole zero cancellation, i.e., for $\hat{b}=b$ it holds that
\begin{align*}
{ G^{-1}_{p_\text{vsi}, \, \omega_\text{sm}}(s) = \dfrac{\tau\, \rho}{\alpha_g (\tau-\rho)} s + k_p + \dfrac{k_i}{s}.}
\end{align*}
In contrast, for $\hat{b}\neq b$ we obtain 
\begin{align*}
 G^{-1}_{p_\text{vsi}, \, \omega_\text{sm}}(s) &=\underbrace{ \left(\dfrac{\tau\, \rho}{\alpha_g (\tau-\rho)} \!+\! \left(\dfrac{1}{b}-\dfrac{1}{\hat{b}}\right)\right) }_{\eqqcolon \xi} s+ k_p + \dfrac{k_i}{s}.
\end{align*}
We are now ready to provide stability conditions for the case $\hat{b}\neq b$. To this end, we define the set
\begin{align}
    \!\mathcal{M} \!\coloneqq\!  \left\{ \begin{bmatrix}\rho\\ \hat{b}\end{bmatrix} \in \mathbb{R}^2 \,\bigg\vert\, \hat{b} < b, \dfrac{\tau\,\alpha_g\,(b-\hat{b})}{\tau\,\hat{b}\,b+\alpha_g\,(b-\hat{b})} \leq \rho < \tau \right\}\!,\!
\end{align} 
and note that the closed-loop transfer function no longer matches the target if $\hat{b}\neq b$, i.e., $G^{\text{cl}}_{\omega_\text{sm}, \, p_\ell}(s) \neq G^{\text{cl}\star}_{\omega_\text{sm}, \, p_\ell}(s)$.

\begin{theorem}
The closed-loop transfer function $G^{\text{cl}}_{\omega_\text{sm}, \, p_\ell}(s)$ is stable if either (i) $\hat{b} \geq b$, or (ii) $(\rho,\hat{b}) \in \mathcal{M}$.
\end{theorem}

\begin{proof}
The gains $k_p$ and $k_i$ defined in \eqref{eq:pid_gains} are positive irrespective of $\hat{b}$. Moreover, for $\hat{b} \geq b$ the coefficient $\xi$ is positive and it directly follows from \cite[Cor. 11]{CS2009} that is $G_{\omega_\text{sm}, \, p_\text{vsi}}(s)$ is positive real. Similarly, if $(\rho,\hat{b}) \in \mathcal{M}$ it can be verified that $\xi$ is positive and $G_{\omega_\text{sm}, \, p_\text{vsi}}(s)$ is positive real. Thus, the sum of the turbine dynamics and $G_{\omega_\text{sm}, \, p_\text{vsi}}(s)$ is obtained as a sum of two passive transfer functions (as depicted in Figure~\ref{fig:inv_interconnection}). As in Proposition~\ref{propo:stablesys}, the closed-loop transfer function ${G}^{\text{cl}}_{\omega_\text{sm}, \, p_\ell}(s)$ is stable as it is a negative feedback interconnection of a strictly passive system and a passive system.
\end{proof}

Moreover, note that $\{\rho \in \mathbb{R} \vert\, \exists\, \,\hat{b}: (\rho,\hat{b}) \in \mathcal{M}\} \subset \mathcal{N} \subset \mathcal{U}$, i.e., if $\hat{b}<b$, the range of effective turbine time constants for which we can guarantee stability has to be further restricted. Therefore, we can also conclude that an overestimate for the susceptance, while degrading the system performance, does not result in an unstable system, while with an underestimate, too small a choice of $\rho$ can result in an unstable system.

Next, we quantify the sensitivity of the closed-loop system to using  $\hat{b}\neq b$ in the controller. To simplify the analysis, we introduce the scaled inverse susceptance mismatch 
\begin{align}\label{eq:beta}
\beta\coloneqq\alpha_g (\tau-\rho)\left(\dfrac{1}{b}-\dfrac{1}{\hat{b}}\right).
\end{align}

\begin{theorem}
\label{th:bounds}
Given a mismatch $\beta$, the $\mathcal H_\infty$-norm of the sensitivity 
\begin{align*}
S(s) \coloneqq \dfrac{\partial}{\partial \beta} 
\dfrac{G^{\text{cl}}_{\omega_\text{sm}, \, p_\ell}(s) - G^{\text{cl}\star}_{\omega_\text{sm}, \, p_\ell}(s)}{G^{\text{cl}\star}_{\omega_\text{sm}, \, p_\ell}(s)}\bigg\vert_{\beta=0}
\end{align*}
of the relative mismatch of the closed-loop transfer function $G^{\text{cl}}_{\omega_\text{sm}, \, p_\ell}(s)$ to a variation in $\beta$ around $0$ can be expressed as 
\begin{align*}
\| S \|_\infty\, =\dfrac{\alpha_g\,(\tau-\rho)}{(\tau^2\rho)}.
\end{align*}
\end{theorem}

\begin{proof}
Using \eqref{eq:beta}, \eqref{eq:overallresp} can be rewritten as
\begin{align*}
    G_{\omega_\text{sm},\, p_\ell}^\text{cl}(s)=\dfrac{-1}{ms+d + \dfrac{\alpha_g [\beta s^2 + (s\,\tau+1)^2]}{(s\,\tau+1)[ \beta s^2 + (s\,\tau +1)(s\,\rho+1)]}}.
\end{align*}
A lengthy calculation reveals the following expression for $S$
\begin{align*}
S(s) = \dfrac{\alpha_g\,(\tau-\rho)\,s^3}{(s\,\rho+1)(s\,\tau +1 )^2},
\end{align*}
and consequently that $S(s)$ has a real pole with multiplicity of two at $-\tfrac{1}{\tau}$ and a real pole multiplicity of one at $-\tfrac{1}{\rho}$. It directly follows that the magnitude $|S(j \omega)|$ is monotonically increasing in $\omega$ and the maximum gain is attained at $\omega \to \infty$ and the Theorem directly follows from
\begin{align*}
\|S\|_\infty\! =\! \lim_{\omega \to \infty} \dfrac{\alpha_g\,(\tau-\rho)\,\omega^3}{\sqrt{(1+\omega^2\,\rho^2)}\,(1+\omega^2\,\tau^2)}\!\! =\! \dfrac{\alpha_g\,(\tau-\rho)}{(\tau^2\rho)}.
\end{align*}
\end{proof}

The key inference drawn from the above theorem (while restricting ourselves to variations around $\beta=0$) is that the {peak} sensitivity $\| S \|_\infty$ of the relative mismatch between $G^{\text{cl}}_{\omega_\text{sm}, \, p_\ell}(s)$ and $G^{\text{cl}\star}_{\omega_\text{sm}, \, p_\ell}(s)$ with respect to the scaled inverse susceptance mismatch $\beta$ is independent of $b$, $\hat{b}$ and  decreases as the tuning parameter $\rho$ increases  (i.e., $\rho \to \tau$). Therefore, increasing $\rho$ decreases the relative mismatch between $G^{\text{cl}}_{\omega_\text{sm}, \, p_\ell}(s)$ and $G^{\text{cl}\star}_{\omega_\text{sm}, \, p_\ell}(s)$. To illustrate how $\beta$ changes with the susceptance mismatch, we let $\hat{b}=c\, b$, thus $\beta\propto \tfrac{c-1}{cb}$, i.e., for a fixed $c$ and effective time constant $\rho$, the  scaled inverse susceptance mismatch $\beta$ and overall mismatch between $G^{\text{cl}}_{\omega_\text{sm}, \, p_\ell}(s)$ and $G^{\text{cl}\star}_{\omega_\text{sm}, \, p_\ell}(s)$ decrease as the true susceptance $b$ increases.

\subsection{Two-bus system case study}
We revisit the two-bus system to illustrate that the impact of incorrect estimates as predicted by Theorem~\ref{th:bounds}. In Figure~\ref{fig:difffreq}, we consider the frequency response to a $1~\mathrm{p.u.}$ load perturbation under two scenarios {\it (i)} $\hat{b}=b$, {\it (ii)} $\hat{b}\neq b$. We note that the mismatch results in a significant variation from the desired trajectory, resulting in a much higher frequency nadir.

\begin{figure}[h]
\centering
\tikzset{every picture/.style={scale=0.99}}%
\input{freq_resp_c.tex}
\caption{The frequency response for different values of $c$, with $\rho=0.7$s.}
\label{fig:difffreq}
\medskip
\centering
\tikzset{every picture/.style={scale=1.0}}%
%
%

\definecolor{mycolor1}{HTML}{BE9063}%
\definecolor{mycolor2}{RGB}{168,50,45}
\definecolor{mycolor3}{HTML}{525B56}
\definecolor{mycolor4}{HTML}{A4978E}%
\definecolor{mycolor5}{rgb}{0.59400,0.18400,0.35600}%
\definecolor{mycolor6}{HTML}{253F5B}%
\definecolor{mycolor7}{HTML}{818A6F}
\definecolor{mycolor8}{HTML}{D59B2D}
\definecolor{mycolor9}{RGB}{31,140,24}
\definecolor{mycolor10}{RGB}{136,176,197}

\begin{tikzpicture}

\begin{axis}[%
width=2.9in,
height=1.2in,
at={(1.739in,0.849in)},
scale only axis,
xmode=log,
xmin=0.001,
xmax=1000,
xminorticks=true,
ymin=260,
ymax=470,
ymajorgrids,
xmajorgrids,
yticklabel style = {font=\footnotesize,xshift=0ex},
xticklabel style = {font=\footnotesize,yshift=0ex},
ylabel style={font=\color{black}},
ylabel={\small Frequency nadir (mHz)},
xlabel style={font=\color{black}},
xlabel={\small True susceptance (p.u.)},
axis background/.style={fill=white},
legend style={legend cell align=left, align=left, draw=black, font=\small, draw=none, legend columns=-1, at={(0.87,0.98)}}
]

\addplot [color=mycolor1, line width=1.75pt]
  table[row sep=crcr]{%
0.001	424.264232792861\\
0.01	424.060818809786\\
0.1	422.242716065556\\
1	414.696563926664\\
10	408.635403188877\\
100	407.612584100144\\
1000	407.503858110404\\
};
\addlegendentry{$0.9$s}

\addplot [color=mycolor2, line width=1.75pt]
  table[row sep=crcr]{%
0.001	424.264213261478\\
0.01	424.058891172688\\
0.1	422.072789148691\\
1	408.739745000999\\
10	380.727516285262\\
100	371.45013184512\\
1000	370.331322668487\\
};
\addlegendentry{$0.7$s}

\addplot [color=mycolor3, line width=1.75pt]
  table[row sep=crcr]{%
0.001	424.264209351358\\
0.01	424.058501855388\\
0.1	422.035494250636\\
1	406.524772456741\\
10	357.487944905099\\
100	331.098666106028\\
1000	327.271118829693\\
};
\addlegendentry{$0.5$s}

\addplot [color=mycolor4, line width=1.75pt]
  table[row sep=crcr]{%
0.001	424.2642076752\\
0.01	424.058334614443\\
0.1	422.019138632622\\
1	405.371370239239\\
10	338.104786185099\\
100	284.577732294178\\
1000	275.360561421446\\
};
\addlegendentry{$0.3$s}

\addplot [color=mycolor1, line width=1.75pt, dashed]
  table[row sep=crcr]{%
0.001	407.491702125468\\
0.01	407.491702125468\\
0.1	407.491702125468\\
1	407.491702125468\\
10	407.491702125468\\
100	407.491702125468\\
1000	407.491702125468\\
};

\addplot [color=mycolor2, line width=1.75pt, dashed]
  table[row sep=crcr]{%
0.001	370.204252215207\\
0.01	370.204252215207\\
0.1	370.204252215207\\
1	370.204252215207\\
10	370.204252215207\\
100	370.204252215207\\
1000	370.204252215207\\
};

\addplot [color=mycolor3, line width=1.75pt, dashed]
  table[row sep=crcr]{%
0.001	326.830556963332\\
0.01	326.830556963332\\
0.1	326.830556963332\\
1	326.830556963332\\
10	326.830556963332\\
100	326.830556963332\\
1000	326.830556963332\\
};

\addplot [color=mycolor4, line width=1.75pt, dashed]
  table[row sep=crcr]{%
0.001	274.299239329325\\
0.01	274.299239329325\\
0.1	274.299239329325\\
1	274.299239329325\\
10	274.299239329325\\
100	274.299239329325\\
1000	274.299239329325\\
};

\end{axis}
\end{tikzpicture}%
\caption{The effect of susceptance mismatch for different values of $\rho$, with $c=1.05$. The dashed lines indicate the frequency nadir for $c=1$.}
\label{fig:difffreqrho}
\medskip
\centering
\tikzset{every picture/.style={scale=1.0}}%
%
%

\definecolor{mycolor1}{HTML}{BE9063}%
\definecolor{mycolor2}{RGB}{168,50,45}
\definecolor{mycolor3}{HTML}{525B56}
\definecolor{mycolor4}{HTML}{A4978E}%
\definecolor{mycolor5}{rgb}{0.59400,0.18400,0.35600}%
\definecolor{mycolor6}{HTML}{253F5B}%
\definecolor{mycolor7}{HTML}{818A6F}
\definecolor{mycolor8}{HTML}{D59B2D}
\definecolor{mycolor9}{RGB}{31,140,24}
\definecolor{mycolor10}{RGB}{136,176,197}

\begin{tikzpicture}

\begin{axis}[%
width=2.9in,
height=1.2in,
at={(1.739in,0.849in)},
scale only axis,
xmode=log,
xmin=0.001,
xmax=1000,
xminorticks=true,
ymin=360,
ymax=440,
ymajorgrids,
xmajorgrids,
yticklabel style = {font=\footnotesize,xshift=0ex},
xticklabel style = {font=\footnotesize,yshift=0ex},
ylabel style={font=\color{black}},
ylabel={\small Frequency nadir (mHz)},
xlabel style={font=\color{black}},
xlabel={\small True susceptance (p.u.)},
axis background/.style={fill=white},
legend style={legend cell align=left, align=left, draw=black, font=\small, draw=none, legend columns=-1, at={(0.95,0.999)}}
]
\addplot [color=mycolor1, line width=1.75pt]
  table[row sep=crcr]{%
0.001	370.204252215207\\
0.01	370.204252215207\\
0.1	370.204252215207\\
1	370.204252215207\\
10	370.204252215207\\
100	370.204252215207\\
1000	370.204252215207\\
};
\addlegendentry{$c=1$}

\addplot [color=mycolor3, line width=1.75pt]
  table[row sep=crcr]{%
0.001	424.177150961324\\
0.01	423.203539025587\\
0.1	415.532000209311\\
1	388.55385949399\\
10	372.734234605693\\
100	370.467833852335\\
1000	370.230718709008\\
};
\addlegendentry{$1.01$}

\addplot [color=mycolor2, line width=1.75pt]
  table[row sep=crcr]{%
0.001	424.264213261478\\
0.01	424.058891172688\\
0.1	422.072789148691\\
1	408.739745000999\\
10	380.727516285262\\
100	371.45013184512\\
1000	370.331322668487\\
};
\addlegendentry{$1.05$}

\addplot [color=mycolor4, line width=1.75pt]
  table[row sep=crcr]{%
0.001	424.284923769715\\
0.01	424.265302970019\\
0.1	424.06972331703\\
1	422.17490585624\\
10	409.286798596004\\
100	381.164883145282\\
1000	371.510993210161\\
};
\addlegendentry{$2$}

\addplot [color=mycolor7, line width=1.75pt]
  table[row sep=crcr]{%
0.001	424.285741551711\\
0.01	424.273476910881\\
0.1	424.151076259716\\
1	422.951187096172\\
10	413.864514185686\\
100	386.013498186105\\
1000	372.266393687757\\
};
\addlegendentry{$5$}

\end{axis}

\end{tikzpicture}%
\caption{The effect of susceptance mismatch with $\rho=0.7$s, $\hat{b}=c\,b$.}
\label{fig:difffreqnadir}
\end{figure}

Next, we analyze the impact of the effective turbine constant $\rho$ on the sensitivity. We recall from Theorem~\ref{th:bounds} that increasing $\rho$ results in a lower sensitivity and decreases the mismatch between the actual and target dynamics. This is illustrated in Figure~\ref{fig:difffreqrho} for varying values of $\rho$.
Finally, we consider the effect of larger mismatch errors for a fixed $\rho$. From Theorem~\ref{th:bounds}, the relative mismatch between $G^{\text{cl}}_{\omega_\text{sm}, \, p_\ell}(s)$ and $G^{\text{cl}\star}_{\omega_\text{sm}, \, p_\ell}(s)$ scales inversely with $b$ around $\beta=0$. However, this trend is observed even for large mismatch values of $\beta$ (i.e., $c>1.05$) as illustrated in Figure~\ref{fig:difffreqnadir}. Finally, the mismatch $c$ has a significant impact on weakly-coupled systems while tightly-coupled systems are less impacted.

\section{Conclusions}
We proposed and analyzed a grid-forming frequency shaping control with a second-order target behaviour for weakly coupled low-inertia power systems. The proposed approach relaxes several existing assumptions in the literature and allows to trade-off IBR peak power injection and SG frequency nadir. Furthermore, we highlighted the significant role of the network parameters and analyzed the effect on the closed-loop system behaviour due to imperfect knowledge of these parameters. While the analytical results are promising, the robustness of the proposed control needs to be validated in high-fidelity simulations. We aim to fully leverage the flexibility offered by the frequency shaping control and different inverter implementations as future work.

\bibliographystyle{IEEEtran}


\end{document}